# MEASURE CONCENTRATION FOR EUCLIDEAN DISTANCE IN THE CASE OF DEPENDENT RANDOM VARIABLES[1]

BY KATALIN MARTON

*Hungarian Academy of Sciences*

Let $q^n$ be a continuous density function in $n$-dimensional Euclidean space. We think of $q^n$ as the density function of some random sequence $X^n$ with values in $\mathbb{R}^n$. For $I \subset [1,n]$, let $X_I$ denote the collection of coordinates $X_i$, $i \in I$, and let $\bar{X}_I$ denote the collection of coordinates $X_i$, $i \notin I$. We denote by $Q_I(x_I|\bar{x}_I)$ the joint conditional density function of $X_I$, given $\bar{X}_I$. We prove measure concentration for $q^n$ in the case when, for an appropriate class of sets $I$, (i) the conditional densities $Q_I(x_I|\bar{x}_I)$, as functions of $x_I$, uniformly satisfy a logarithmic Sobolev inequality and (ii) these conditional densities also satisfy a contractivity condition related to Dobrushin and Shlosman's strong mixing condition.

**1. Introduction.** Let us consider the absolutely continuous probability measures in $n$-dimensional Euclidean space $\mathbb{R}^n$. With some abuse of notation, we use the same letter to denote a probability measure and its density function.

We say that a measure $q^n$ on $\mathbb{R}^n$ has the measure concentration property (with respect to the Euclidean distance) if

$$(1.1) \quad d(A,B) \leq c \cdot \left[ \sqrt{\log \frac{1}{q^n(A)}} + \sqrt{\log \frac{1}{q^n(B)}} \right] \quad \text{for any sets } A, B \subset \mathbb{R}^n,$$

where $d(A,B)$ denotes the Euclidean distance of the sets $A$ and $B$. (We consider only measurable sets. This definition is equivalent to the more familiar definition that involves the probabilities of a set $A$ and its $\varepsilon$ neighborhood.)

Measure concentration is an important property, since it implies sub-Gaussian behavior of the Laplace transforms of Lipschitz functions and

Received July 2002; revised December 2003.
[1]Supported in part by Hungarian Academy of Sciences Grants OTKA T 26041 and T 32323.

thereby is an important tool for proving strong forms of the law of large numbers.

Measure concentration for $q^n$ follows from the validity of a logarithmic Sobolev inequality for $q^n$ by a recent theorem of Otto and Villani (2000). However, in this paper we prove measure concentration in some cases when a logarithmic Sobolev inequality probably cannot be proved by an easy extension of the existing methods.

Consider the following the distance between probability measures in $\mathbb{R}^n$:

$$W(p^n, q^n) = \inf_\pi [E_\pi (Y^n - X^n)^2]^{1/2},$$

where $Y^n$ and $X^n$ are random variables distributed according to the law $p^n$ and $q^n$, respectively, and the infimum is taken over all distributions $\pi$ on $\mathbb{R}^n \times \mathbb{R}^n$ that have $p^n$ and $q^n$ as marginals. This is one of the transportation cost related distances between measures, often called the Wasserstein distance (based on the squared Euclidean distance).

Let us denote by $D(p^n \| q^n)$ the informational divergence of the probability distribution $p^n$ with respect to $q^n$ as

$$D(p^n \| q^n) = \int_{\mathbb{R}^n} \log \frac{dp^n}{dq^n} \, dp^n$$

if $p^n$ is absolutely continuous with respect to $q^n$, and $\infty$ otherwise. (By $dp^n/dq^n$ we denote the Radon–Nikodym derivative.)

By a simple argument [first used by Marton (1986, 1996) for Hamming distance and then by Talagrand (1996) for Euclidean distance] it can be shown that if $q^n$ satisfies, for some $\rho > 0$, the "distance-divergence" inequality

$$(1.2) \quad W(p^n, q^n) \leq \sqrt{\frac{2D(p^n \| q^n)}{\rho}} \quad \text{for any probability measure } p^n \text{ on } \mathbb{R}^n,$$

then it satisfies the measure concentration inequality (1.1) as well (with $c = \sqrt{2/\rho}$). Indeed, assume (1.2) and let $A, B \subset \mathbb{R}^n$ be measurable sets in $\mathbb{R}^n$. Denote by $p^n$ and $r^n$ the restriction of $q^n$ to $A$ and $B$, respectively:

$$p^n(C) = \frac{q^n(C \cap A)}{q^n(A)}, \qquad r^n(C) = \frac{q^n(C \cap B)}{q^n(B)}.$$

Since $p^n$ and $r^n$ are supported by $A$ and $B$, respectively, we have, using also the triangle inequality for $W$,

$$d(A, B) \leq W(p^n, r^n) \leq W(p^n, q^n) + W(r^n, q^n)$$

$$\leq \sqrt{\frac{2D(p^n \| q^n)}{\rho}} + \sqrt{\frac{2D(r^n \| q^n)}{\rho}}.$$

Since

$$D(p^n \| q^n) = \log \frac{1}{q^n(A)} \quad \text{and} \quad D(r^n \| q^n) = \log \frac{1}{q^n(B)},$$

(1.1) follows.

Therefore, our aim is to find possibly general sufficient conditions for a measure $q^n$ to satisfy a distance-divergence inequality (1.2).

The inequality (1.2) was first proved by Talagrand (1996) for the case when $q^n$ is a Gaussian product measure, and his proof for dimension 1 easily generalizes to the case when $q^1$ is uniformly log-concave. [Second derivative of $-\log q^1(x)$ bounded from below.] With more effort, using recent results on the solution of the Monge–Kantorovich problem for the Wasserstein distance in the Euclidean space [McCann (1995)], this generalization can be carried out for the multidimensional case as well. For an alternative proof, see Bobkov and Götze (1999). Otto and Villani (2000) proved that the distance-divergence inequality (1.2) follows from $q^n$ satisfying a logarithmic Sobolev inequality, that is, from

$$D(p^n\|q^n) \leq \frac{1}{2\rho} \int_{\mathbb{R}^n} \left|\nabla \log \frac{p^n}{q^n}\right|^2 dq^n$$

holding for any density function $p^n$ on $\mathbb{R}^n$ such that $p^n/q^n$ is smooth enough. A simple sufficient condition for $q^n$ satisfying a logarithmic Sobolev inequality is that $q^n$ be a bounded perturbation of a uniformly log-concave function. (See later.)

Much effort has been spent to find sufficient conditions for the logarithmic Sobolev inequality in terms of the conditional density functions

$$Q_i(\cdot|\bar{x}_i) = \text{dist}_{q^n}(X_i|\bar{X}_i = \bar{x}_i)$$

of $q^n$. [Here $\bar{X}_i$ denotes $(X_k : k \neq i)$.] This problem is not yet satisfactorily solved. In the cases considered, $q^n$ is a Gibbs state (with unbounded spins) over a region $\Lambda$ of the $d$-dimensional integer lattice and $q^n$ corresponds to a pair interaction with bounded range. [See Yoshida (1999a, b), Helffer (1999), Bodineau and Helffer (1999) and Ledoux (1999).]

In this paper we use a different approach: To prove distance-divergence inequality for $q^n$, we use the one-dimensional distance-divergence inequality for the conditional distributions $Q_i(\cdot|\bar{x}_i)$.

NOTATION. The integers $i$, $1 \leq i \leq n$, are called sites and $[1,n]$ is the set of sites. Let $\mathcal{I}$ be a family of sets $I \subset [1,n]$, called patches. Each $I \in \mathcal{I}$ has a multiplicity $\geq 1$ and the number of patches counted with multiplicities is denoted by $N$. A patch consisting of one element $i$ is denoted by $i$: for $x^n \in \mathbb{R}^n$ and $J \subset [1,n]$, $x_J = (x_i : i \in J)$ and $\bar{x}_J = (x_i : i \notin J)$; for $a^n \in \mathbb{R}^n$ and $J \subset [1,n]$, $|a_J|^2 = \sum_{i \in J} a_i^2$.

Let $q^n$ denote the density of an absolutely continuous probability measure on $\mathbb{R}^n$ and let $X^n$ denote a random sequence in $\mathbb{R}^n$, dist $X^n = q^n$. Conditional density functions consistent with $q^n$ are expressed as $Q_I(\cdot|\bar{x}_I) = \text{dist}(X_I|\bar{X}_I = \bar{x}_I)$ ($I \in \mathcal{I}$), whereas $\bar{q}_I(\bar{x}_I)$ ($I \in \mathcal{I}$) denotes the density function of $\bar{X}_I$. The density of a probability distribution on $\mathbb{R}^n$ is denoted by $p^n$ and $Y^n$ denotes a random sequence with dist $Y^n = p^n$. Conditional density functions consistent with $p^n$ are expressed $p_I(\cdot|\bar{x}_I) = \mathrm{dist}(Y_I|\bar{Y}_I = \bar{x}_I)$ and the density function of $\bar{Y}_I$ is denoted $\bar{p}_I(\bar{y}_I)$ ($I \in \mathcal{I}$).

Let $\mathcal{W}_I(\rho)$ represent the set of all probability distributions $Q_I$ on $\mathbb{R}^I$ that satisfy, for every distribution $p_I$ on $\mathbb{R}^I$, the distance-divergence inequality

$$W(p_I, Q_I) \leq \sqrt{\frac{2D(p_I \| Q_I)}{\rho}}.$$

In the simplest case, $\mathcal{I}$ consists of the one-element sets of $[1, n]$. Alternatively, let $n$ be the cardinality of a (large) box $\Lambda$ in the $d$-dimensional lattice $\mathbb{Z}^d$, let $V \subset \mathbb{Z}^d$ be a (relatively small) set and let $\mathcal{I}$ consist of the intersections of the translates of $V$ with $\Lambda$. The multiplicity of such a set $I$ can be taken as the number of different translates of $V$ whose intersection with $\Lambda$ is $I$. Every site is covered then by $|V|$ patches, where $|V|$ is the cardinality of $V$.

Theorem 1 presents a sufficient condition for a distance-divergence inequality of type (1.2) in terms of the conditional distributions $Q_I(\cdot|\bar{x}_I)$ ($I \in \mathcal{I}$). The reason we want a condition in terms of the conditional distributions $Q_I(\cdot|\bar{x}_I)$ is that in statistical physics the model is often defined in such a manner that it gives direct information on these conditional distributions. For example, $q^n$ may be the conditional distribution of a Gibbs random field over a domain in a multi-dimensional lattice with fixed boundary condition.

The conditions of the theorem require that the individual conditional distributions $Q_I(\cdot|\bar{x}_I)$ behave nicely, and we also need the following assumption on the ensemble of the conditional distributions $Q_I(\cdot|\bar{x}_I)$ ($I \in \mathcal{I}$).

DEFINITION 1 (Contractivity condition). Let $\mathcal{I}$ be such that every site $i$ is covered by at least $t \geq 1$ patches $I$. The system of conditional distributions $Q_I(\cdot|\bar{x}_I)$ ($I \in \mathcal{I}$) is $(1-\delta)$-contractive ($\delta > 0$) if for any pair of sequences $(y^n, z^n) \in \mathbb{R}^n \times \mathbb{R}^n$,

$$(1.3) \qquad \sum_{I \in \mathcal{I}} W^2(Q_I(\cdot|\bar{y}_I), Q_I(\cdot|\bar{z}_I)) \leq t \cdot (1-\delta)|y^n - z^n|^2.$$

For clarity, we formulate the contractivity condition for the special case when $\mathcal{I}$ is the family of one-element patches.

CONTRACTIVITY CONDITION FOR ONE-ELEMENT PATCHES. We say that the system of conditional distributions $Q_i(\cdot|\bar{x}_i)$ is $(1-\delta)$-contractive ($\delta > 0$)

if for any pair of sequences $(y^n, z^n) \in \mathbb{R}^n \times \mathbb{R}^n$,

$$(1.3') \qquad \sum_{i=1}^{n} W^2(Q_i(\cdot|\bar{y}_i), Q_i(\cdot|\bar{z}_i)) \leq (1-\delta)|y^n - z^n|^2.$$

The contractivity condition is related to Dobrushin and Shlosman's strong mixing condition. Indeed, it is obviously implied by the following condition:

DEFINITION 2 (Dobrushin–Shlosman-type contractivity condition). Let us assume again that every site $i$ is covered by at least $t \geq 1$ patches $I$. We say that the system of conditional distributions $Q_I(\cdot|\bar{x}_I)$ ($I \in \mathcal{I}$) satisfies a Dobrushin–Shlosman-type contractivity condition if for every $I \in \mathcal{I}$ and $k \notin I$, and for every $\bar{y}_I, \bar{x}_I$ differing only at site $k$,

$$(1.4) \qquad W(Q_I(\cdot|\bar{x}_I), Q_I(\cdot|\bar{y}_I)) \leq \alpha_{k,I}|y_k - x_k|,$$

and for the matrix $A = (\alpha_{k,I})$ ($\alpha_{k,I} = 0$ for $k \in I$ by definition)

$$(1.3'') \qquad \|A\|^2 \leq (1-\delta) \cdot t.$$

Here $\|A\|$ denotes the norm of $A$ considered as an $L_2 \mapsto L_2$ operator and $\delta$ is a positive constant.

To see that Definition 2 is stronger than Definition 1, note that, by the triangle inequality, (1.4) implies

$$W(Q_I(\cdot|\bar{y}_I), Q_I(\cdot|\bar{z}_I)) \leq \sum_{k \notin I} \alpha_{k,I}|z_k - y_k| \qquad \text{for all } I, y^n, z^n,$$

so, by the definition of $\|A\|$,

$$\sum_{I \in \mathcal{I}} W^2(Q_I(\cdot|\bar{y}_I), Q_I(\cdot|\bar{z}_I))$$

$$\leq \sum_{I \in \mathcal{I}} \left(\sum_k \alpha_{k,I}|z_k - y_k|\right)^2 \leq \|A\|^2 |y^n - z^n|^2 \leq t \cdot (1-\delta)|y^n - z^n|^2.$$

The stronger version of the contractivity condition given in Definition 2, when considered for one-element patches, is an analog of Dobrushin's (1970) uniqueness condition. In the general case, it is an analog of Dobrushin and Shlosman's (1985a) uniqueness condition (CV). Note, however, that we use a variant of the Wasserstein distance that minimizes the expected squared distance, whereas in condition (CV) in Dobrushin and Shlosman's (1985a) a form of the Wasserstein distance is used that minimizes the expected distance without squaring. Moreover, we require (1.4) to hold for all patches $I \in \mathcal{I}$, which means, for the second example considered above, all the intersections of the translates of a given set with another given set. This is reminiscent of condition (CC) in Yoshida (1999b), who formulated a set of mixing conditions, one of which is (CC), that can be considered the analogs of Dobrushin and Shlosman's (1985b, 1987) strong mixing conditions. It is not completely clear how the contractivity condition used in this paper (Definition 1) is related to the set of mixing conditions in Yoshida (1999b). However, we think that the conditions in Yoshida (1999b) should not be considered final and standard yet, since their equivalence among each other and with the logarithmic Sobolev inequality is only proved for ferromagnetic interactions and in the case of superquadratic growth of the single spin phase. We think that the contractivity condition is understandable in itself, and we do not need an analysis of its analogy with the Dobrushin–Shlosman conditions, or the conditions in Yoshida (1999b). We note, however, that we assume nothing that would correspond to the boundedness of the "range of interaction."

In this paper we use only Definition 1. Definition 2 is stated here only to explain the relationship with previously existing concepts.

THEOREM 1. *Let $\mathcal{I}$ be such that every site $i$ is covered by at least $t \geq 1$ and by at most $v$ patches. Assume that for all $I \in \mathcal{I}$, $Q_I(x_I|\bar{x}_I)$, as a function of $n$ variables, is continuous. Assume further that for every $I \in \mathcal{I}$ and every $\bar{x}_I$,*

$$\text{(1.5)} \qquad Q_I(\cdot|\bar{x}_I) \in \mathcal{W}_I(\rho).$$

*Finally, assume that the system of conditional distributions $Q_I(\cdot|\bar{x}_I)$ $(I \in \mathcal{I})$ is $(1-\delta)$-contractive $(\delta > 0)$. Then for any distribution $p^n$ on $\mathbb{R}^n$,*

$$\text{(1.6)} \qquad W(p^n, q^n) \leq C \cdot \sqrt{\frac{v}{t} \cdot \frac{1}{\delta} \cdot \frac{2}{\rho} \cdot D(p^n \| q^n)},$$

*where $C$ is a numerical constant.*

Formula (1.6) simplifies if $v = t$, as in the above examples.

The conditions of Theorem 1 are quite abstract, so we are going to formulate a special case where the conditions can be verified.

Write the density function $q^n$ in the form

$$\text{(1.7)} \qquad q^n(x^n) = \frac{1}{Z} \cdot \exp(-\Phi(x^n)),$$

where $Z$ is a normalizing constant. Then the conditional density functions are of the form

$$\text{(1.8)} \qquad Q_I(x_I|\bar{x}_I) = \frac{1}{Z(\bar{x}_I)} \cdot \exp(-\Phi(x^n)),$$

where $Z(\bar{x}_I)$ is the normalizing factor.

When does the ensemble of the conditional densities $Q_I(x_I|\bar{x}_I)$ satisfy the conditions of Theorem 1?

It is natural to try to use the recent result by Otto and Villani that deduces the distance-divergence inequality for some probability measure $q$ on $\mathbb{R}^k$, that is, the relationship

$$W(p,q) \leq \sqrt{\frac{2D(p\|q)}{\rho}} \qquad \text{for all } p,$$

from $q$ satisfying a logarithmic Sobolev inequality.

DEFINITION 3. The density function $q$ on $\mathbb{R}^k$ satisfies a logarithmic Sobolev inequality with constant $\rho$ if for any density function $p$ on $\mathbb{R}^k$, such that $p(x^k)/q(x^k)$ is sufficiently smooth,

$$D(p\|q) \leq \frac{1}{2\rho} \int_{\mathbb{R}^k} \left|\nabla \log \frac{dp}{dq}\right|^2 dp.$$

The following sufficient condition follows from the Bakry–Emery (1985) criterion, supplemented by a perturbation result from Holley and Stroock (1987):

PROPOSITION 1. *Let $q(x^k)$ be a density function of the form $\exp[-V(x^k)]$ and let $V$ be strictly convex at $\infty$, that is, $V(x^k) = U(x^k) + K(x^k)$, where $K(x^k)$ is bounded, and the Hessian*

$$D(x) = (\partial_{ij} U(x))$$

*satisfies*

$$D(x) \geq c \cdot I$$

*for some $c > 0$ (where $I$ is the identity matrix). Then $q$ satisfies a logarithmic Sobolev inequality with constant $\rho$, depending only on $c$ and $\|K\|_\infty$:*

$$\rho \geq c \cdot \exp(-4\|K\|_\infty).$$

Note that, on the real line, a necessary and sufficient condition for a density function to satisfy a logarithmic Sobolev inequality was established by Bobkov and Götze. From this result, Gentil (2001) derived the logarithmic Sobolev inequality for a class of density functions, different from the above class. We do not cite this theorem.

THEOREM OF OTTO AND VILLANI (2000) [see Bobkov, Gentil and Ledoux (2001) also]. *If the density function $q(x^k)$ ($x^k \in \mathbb{R}^k$) satisfies a logarithmic Sobolev inequality with constant $\rho$, then $q \in \mathcal{W}_{[1,k]}(\rho)$.*

[The theorem as stated here was proved by Bobkov, Gentil and Ledoux (2001); its original version from Otto and Villani (2000) contained some minor additional condition.]

To formulate a sufficient condition for the contractivity condition we introduce some notation.

Let $\mathcal{I}$ be a family of patches as in Theorem 1. Consider distribution (1.7) and assume that $\Phi$ is twicely continuously differentiable. For a fixed sequence $y^n \in \mathbb{R}^n$ and a vector $\eta = (\eta_I, I \in \mathcal{I})$, where $\eta_I \in \mathbb{R}^{|I|}$, we define a matrix $B = B(\eta, y^n)$. The rows of $B$ are indexed by pairs $(I, i)$, $(i \in I \in \mathcal{I})$, while its columns are indexed by $k$ $(1 \le k \le n)$,

$$B = B(\eta, y^n) = (\beta_{(I,i),k}(\eta, y^n)),$$

where

$$\beta_{(I,i),k}(\eta, y^n) = \begin{cases} \partial_{ik}\Phi(\eta_I, \bar{y}_I), & i \in I,\ k \notin I, \\ 0, & i,\ k \in I. \end{cases}$$

For the case of one-element patches, the definition of $B(\eta, y^n) = B(\eta^n, y^n)$ becomes quite simple:

$$B = B(\eta^n, y^n) = (\beta_{i,k}(\eta^n, y^n)),$$
$$\beta_{i,k}(\eta^n, y^n) = \begin{cases} \partial_{ik}\Phi(\eta_i, \bar{y}_i), & k \ne i, \\ 0, & k = i. \end{cases}$$

Note that if $\Phi$ has the form

$$\Phi(x^n) = \sum_{i=1}^n V(x_i) + \sum_{i \ne k} b_{i,k} x_i x_k,$$

then $B$ does not depend on $\eta$ and $y^n$. For example, in the case of one-element patches we have

$$B = (b_{i,k}).$$

THEOREM 2. *Let $\mathcal{I}$ be a family of patches as in Theorem 1 and assume that $\Phi$ is twice continuously differentiable. Assume furthermore that the conditional densities (1.8), as functions of $x_I$, satisfy a logarithmic Sobolev inequality with the same $\rho$ (independently of $I$ and $\bar{x}_I$). If*

$$(1.9) \qquad \sup_{\eta, y^n} \left\| \frac{1}{\rho} \cdot B(\eta, y^n) \right\|^2 \le t \cdot (1 - \delta),$$

*then Theorem 1 holds.*

In view of the Otto–Villani theorem, the question arises whether the conditions of Theorem 2 might imply a logarithmic Sobolev inequality.

This is not to be expected with the existing proofs of logarithmic Sobolev inequality for Gibbs fields. Indeed, consider a Gibbs field over a cube $\Lambda \subset \mathbb{Z}^d$, with single spin space $\mathbb{R}$ and potential

$$(1.10) \qquad \Phi(x^n) = \sum_{i \in \Lambda} V(x_i) + \sum_{i,j \in \Lambda} b_{i,j} x_i x_j + \sum_{\substack{i \in \Lambda \\ j \in \mathbb{Z}^d - \Lambda}} b_{i,j} x_i \omega_j,$$

where $\{\omega_j : j \in \mathbb{Z}^d - \Lambda\}$ is the configuration outside $\Lambda$. We assume that $V(x)$ is convex at $\infty$, that is, $V(x) = U(x) + K(x)$, where $U''(x) \geq c > 0$ and $K(x)$ is bounded. If $b_{i,j}$ does not go to 0 exponentially fast with $|i - j| \to \infty$, then the proofs of Yoshida and Bodineau–Helffer break down, whereas it is still possible that condition (1.9) of Theorem 2 holds. If $b_{i,j} = J > 0$ for $i$ and $j$ nearest neighbors, and $b_{i,j} = 0$ otherwise, then Yoshida's proof requires superquadratic growth for the single spin phase $V(x)$ at $\infty$; and for the Bodineau–Helffer proof to work, $2dJ$ must not approach $\rho$, whereas for condition (1.9) of Theorem 2 to hold, it is sufficient that $2dJ < \rho$.

On the other hand, Ledoux's (1999) proof of his Proposition 2.3 does apply for nearest neighbor interactions with interaction coefficient $J > 0$ satisfying $2dJ < c \cdot \exp(-4\|K\|_\infty)$ and proves the correlation bound (DS3) of Yoshida (1999b). By the results of Yoshida (1999b), this is equivalent to the logarithmic Sobolev inequality, provided $V(x)$ grows superquadratically at $\infty$.

However, even if the interaction coefficients and the single spin phase are such that a logarithmic Sobolev inequality holds, that does not yield a simple explicit bound for the logarithmic Sobolev constant nor for the coefficient in the distance-divergence inequality. On the other hand, Theorem 2 implies the following corollary for potential (1.10):

COROLLARY 1. *If for the potential* (1.10), *the density function* const. $\exp(-V(x))$ *satisfies a logarithmic Sobolev inequality with constant $\rho$ and for the (infinite) matrix $B = (b_{i,j})$,*

$$\|B\| < \rho,$$

*then*

$$W^2(p^n, q^n) \leq C^2 \cdot \frac{1}{1 - \|B/\rho\|^2} \cdot \frac{2D(p^n \| q^n)}{\rho}.$$

REMARK 1. It is not known whether the distance-divergence inequality (1.2) implies a logarithmic Sobolev inequality (possibly with a different constant [Villani (2003)] this would be a converse to the Otto–Villani theorem). Thus it would be very interesting to prove or disprove that a logarithmic Sobolev inequality holds under the conditions of Theorem 2, with a constant depending on $\rho$ and $\delta$.

The proof of Theorem 1 is based on a Markov chain (sometimes called the Gibbs sampler), which realizes a discrete time interpolation between $p^n$ and the Markov chain's invariant measure $q^n$. The contractivity condition allows us to prove that this Markov chain converges to $q^n$ exponentially with respect to the Wasserstein distance. (See the end of Section 2.) Before coming to this step, we prove a bound for $W^2(p^n, p^n \Gamma^M)$ (where $p^n \Gamma^M$ is the distribution of the $M$th term of the Markov chain) in terms of $D(p^n \| q^n)$.

**2. Some Markov kernels and probability distributions on $\mathbb{R}^n$.** The following Markov kernels, which are associated with the conditional density functions $Q_I(\cdot | \bar{x}_I)$ ($I \in \mathcal{I}$), are instrumental in our forthcoming constructions.

For $I \in \mathcal{I}$ define the Markov kernel (i.e., the conditional distribution) $\Gamma_I(dz^n | y^n)$ as follows. The projection of $\Gamma_I(\cdot | y^n)$ on the coordinates outside $I$ is defined as

$$\Gamma_I(\{\bar{y}_I\} | y^n) = 1.$$

The projection of $\Gamma_I(\cdot | y^n)$ on the coordinates in $I$ is given by the conditional density $Q_I(\cdot | \bar{y}_I)$:

$$\Gamma_I(dz_I | y^n) = Q_I(z_I | \bar{y}_I) \, dz_I.$$

We define the Markov kernel $\Gamma(dz^n | y^n)$ as a mixture of the Markov kernels $\Gamma_I(\cdot | y^n)$:

$$\Gamma(dz^n | y^n) = \frac{1}{N} \sum_{I \in \mathcal{I}} \Gamma_I(dz^n | y^n).$$

Finally, for an integer $M \geq 0$, we denote by $\Gamma^M(dz^n | y^n)$ the $M$th operator power of $\Gamma(dz^n | y^n)$:

$$\Gamma^M(dz^n | y^n)$$
$$= \int\int \cdots \int \Gamma(dz^n | z^n(M-1)) \Gamma(dz^n(M-1) | z^n(M-2)) \cdots \Gamma(dz^n(1) | y^n).$$

Equivalently,
$$\Gamma^M(dz^n | y^n)$$
$$= \frac{1}{N^M} \sum_{I_1, I_2, \ldots, I_M \in \mathcal{I}} \int\int \cdots \int \Gamma_{I_M}(dz^n | z^n(M-1))$$
$$\times \Gamma_{I_{M-1}}(dz^n(M-1) | z^n(M-2)) \cdots$$
$$\times \Gamma_{I_1}(dz^n(1) | y^n).$$

Now, at the risk of redundancy, we give a somewhat lengthy description of the Markov kernel $\Gamma^M(dz^n|y^n)$, along with some associated density functions, since it is this description that we use in the sequel. We keep $M$ fixed.

Let us fix a sequence of patches

(2.1) $$(I_1, I_2, \ldots, I_M).$$

Also, fix a density function $p^n = \operatorname{dist} Y^n$. We define successively the density functions

(2.2) $$r^n(0) = p^n, \qquad r^n(l) = r^n(l-1)\Gamma_{I_l}, \qquad l = 1, 2, \ldots, M.$$

We think of the density functions $r^n(l)$ as being conditional density functions of random sequences $Z^n(l)$, $l = 1, 2, \ldots, M$, given that in a random independent $M$-wise selection from the set $\mathcal{I}$, we have drawn $I_1, I_2, \ldots, I_M$:

$$r^n(l) = \operatorname{dist}(Z^n(l)|I_1, I_2, \ldots, I_M).$$

It follows from (2.2) that $r^n(l)$ does not depend on $I_{l+1}, \ldots, I_M$, that is,

$$r^n(l) = \operatorname{dist}(Z^n(l)|I_1, I_2, \ldots, I_l)$$

for every $l$.

We also define a joint conditional distribution for $(Z^n(0) = Y^n, Z^n(1), \ldots, Z^n(M))$, given by (2.1). First we define, for every $l$,

$$\operatorname{dist}(Z^n(l-1), Z^n(l)|I_1, I_2, \ldots, I_M) = \operatorname{dist}(Z^n(l-1), Z^n(l)|I_1, I_2, \ldots, I_l)$$

in such a way that

$$\operatorname{dist}(\bar{Z}_{I_l}(l)|Z^n(l-1) = z^n(l-1), I_1, I_2, \ldots, I_l)$$

is concentrated on $\{\bar{z}_{I_l}(l-1)\}$ and, moreover, that

$$\operatorname{dist}(Z_{I_l}(l-1), Z_{I_l}(l)|\bar{Z}_{I_l}(l-1) = \bar{z}_{I_l}(l-1), I_1, I_2, \ldots, I_l)$$

minimizes, for each value of $\bar{z}_{I_l}(l-1)$, the expected conditional quadratic distance

$$E\{|Z_{I_l}(l) - Z_{I_l}(l-1)|^2|\bar{z}_{I_l}(l-1), I_1, I_2, \ldots, I_l\}.$$

At this point we use the condition (1.5) in Theorem 1 to infer that this minimization yields

(2.3) $$E\{|Z_{I_l}(l) - Z_{I_l}(l-1)|^2|\bar{z}_{I_l}(l-1), I_1, I_2, \ldots, I_l\}$$
$$\leq \frac{2}{\rho} \cdot D(r_{I_l}(l-1)(\cdot|\bar{z}_{I_l}(l-1))\|Q_{I_l}(\cdot|\bar{z}_{I_l}(l-1)))$$

for all $\bar{z}_{I_l}(l-1)$.

Finally, we define

$$\operatorname{dist}(Z^n(0), Z^n(1), \ldots, Z^n(M)|I_1, I_2, \ldots, I_M)$$

so that, for $(I_1, I_2, \ldots, I_M)$ fixed, $(Z^n(0), Z^n(1), \ldots, Z^n(M))$ is a Markov chain.

Note that although $r^n(l) = \text{dist}(Z^n(l)|I_1, I_2, \ldots, I_l) = r^n(l-1)\Gamma_{I_l}$,

$$\text{dist}(Z^n(l)|Z^n(l-1), I_1, \ldots, I_l) \neq \Gamma_{I_l}.$$

Taking average with respect to $I_1, I_2, \ldots, I_M$, we get the (unconditional) joint distribution

$$\text{dist}(Z^n(0), Z^n(1), \ldots, Z^n(M)).$$

It is easy to see that

$$\text{dist}\, Z^n(l) = p^n \Gamma^l \qquad \text{for } l = 0, 1, \ldots, M.$$

We use the notation $Y^n$ for $Z^n(0)$ and use $Z^n$ for $Z^n(M)$.

It is important in the sequel that the Markov kernels $\Gamma_I$, $\Gamma$ and $\Gamma^M$ all have $q^n$ as invariant measure.

Note that we could (and do, in fact) consider the infinite Markov chain with marginal distributions $p^n \Gamma^l$, $0 \leq l < \infty$, as well. This infinite Markov chain is a variant of the so-called Gibbs sampler, which is well known in Markov chain simulation.

In Section 4 we prove that $\Gamma$ is a contraction with respect to the Wasserstein distance, which implies that $p^n \Gamma^m \to q^n$ as $m \to \infty$, exponentially fast:

PROPOSITION 2. *Assume that the conditional distribution functions $Q_I(\cdot|\bar{x}_I)$, $I \in \mathcal{I}$, satisfy the contractivity condition (Definition 1) and that every site is covered by at least $t$ patches. Let $p^n = \text{dist}\, Y^n$ and $r^n = \text{dist}\, U^n$ be two density functions on $\mathbb{R}^n$. Then*

$$W^2(p^n\Gamma, r^n\Gamma) \leq \left(1 - \frac{t\delta}{N}\right) \cdot W^2(p^n, r^n).$$

COROLLARY 2. *Under the conditions of Proposition 3,*

$$W^2(p^n\Gamma^m, q^n) \leq \left(1 - \frac{t\delta}{N}\right)^m \cdot W^2(p^n, q^n)$$

*for any integer $m \geq 0$.*

**3. An auxiliary theorem.** A basic tool in the proof of Theorem 1 is the following theorem, which gives a bound for $W(p^n, p^n\Gamma^M)$ in terms of the informational divergence $D(p^n \| q^n)$. We hope that it turns out to be interesting in its own right. For this auxiliary theorem we do not use the contractivity condition.

AUXILIARY THEOREM. *Assume that for every $I \in \mathcal{I}$ and every $\bar{x}_I$, the conditional density function $Q_I(\cdot|\bar{x}_I)$ satisfies condition (1.5), and that each site is covered by at most $v$ patches. Then for any density function $p^n$ on $\mathbb{R}^n$ and for the Markov kernel $\Gamma$,*

$$W^2(p^n, p^n \Gamma^M) \leq \frac{M}{N} \cdot v \cdot \frac{2}{\rho} \cdot D(p^n \| q^n) \qquad \text{for any } M.$$

REMARK 2. For the joint distribution $\text{dist}(Y^n, Z^n)$, with marginals $p^n$ and $p^n \Gamma^M$, yielding $W(p^n, p^n \Gamma^M)$,

$$\text{dist}(Z^n | Y^n) \neq \Gamma^M,$$

in general.

By the construction of the Markov chain $(Y^n = Z^n(0), Z^n(1), \ldots, Z^n(M) = Z^n)$ we have

$$\text{dist } Z^n = p^n \Gamma^M$$

and we use the joint distribution of the Markov chain to estimate $W^2(p^n, p^n \Gamma^M)$. Clearly,

$$W^2(p^n, p^n \Gamma^M) \leq E|Y^n - Z^n|^2.$$

First we prove the following lemma.

LEMMA 1. *We have*

$$E|Y^n - Z^n|^2 \leq \frac{M}{N} \cdot v \cdot \sum_{l=1}^{M} E|Z_{I_l}(l) - Z_{I_l}(l-1)|^2.$$

*(Note that in this formula the subscripts $I_l$ are random and the expected value takes an average with respect to them, too.)*

PROOF OF LEMMA 1. For a realization of the sequence of patches, say

(3.1) $$(I_1, I_2, \ldots, I_M),$$

we denote by $\sigma$ the listing of the sites in the patches (3.1),

$$\sigma = (i_1, i_2, \ldots, i_m, \ldots, i_L), \qquad L = \sum_{l=1}^{M} |I_l|,$$

where $|I_l|$ denotes the cardinality of $I_l$ and $i_m = i \in [1, n]$ if

$$m = \sum_{j=1}^{l-1} |I_j| + r, \qquad 0 < r \leq |I_l|,$$

and the $r$th site in the patch $I_l$ is just $i$. Let $\nu_i$ denote the frequency of $i$ in $\sigma$ and let $\mu_{i,m}$ denote the frequency of $i$ in $(i_1, i_2, \ldots, i_m)$.

Write
$$\zeta_{ij}^2 = |Z_i(l) - Z_i(l-1)|^2, \qquad 1 \leq i \leq n, 1 \leq j \leq \nu_i,$$
if $I_l$ is the $j$th patch in (3.1) that contains the site $i$.

It follows from the triangle and the Cauchy–Schwarz inequalities that
$$(3.2) \qquad E|Y^n - Z^n|^2 \leq \sum_{i=1}^{n} \sum_{k} \Pr\{\nu_i = k\} \cdot k \cdot \sum_{j=1}^{k} E\{\zeta_{ij}^2 | \nu_i = k\}.$$

For $j \leq k$ we have
$$E\{\zeta_{ij}^2 | \nu_i = k\} = E\{\zeta_{ij}^2 | \nu_i \geq j, \nu_i = k\},$$
but $\zeta_{ij}^2$ is conditionally independent of $\nu_i$ under the condition $\{\nu_i \geq j\}$. It follows that for $j \leq k$,
$$E\{\zeta_{ij}^2 | \nu_i = k\} = E\{\zeta_{ij}^2 | \nu_i \geq j\}.$$

Thus (3.2) can be continued to
$$E|Y^n - Z^n|^2 \leq \sum_{i=1}^{n} \sum_{j \geq 1} E\{\zeta_{ij}^2 | \nu_i \geq j\} \cdot \sum_{k \geq j} \Pr\{\nu_i = k\} \cdot k.$$

Furthermore, for any $i, j$,
$$\sum_{k \geq j} \Pr\{\nu_i = k\} \cdot k \leq E\nu_i \leq \frac{M}{N} \cdot v,$$
whence
$$(3.3) \qquad E|Y^n - Z^n|^2 \leq \frac{M}{N} \cdot v \cdot \sum_{i=1}^{n} \sum_{j \geq 1} E\{\zeta_{ij}^2 | \nu_i \geq j\}.$$

Put
$$\eta_m^2 = \zeta_{ij}^2, \qquad m = 1, 2, \ldots, L,$$
where $(i, j)$ and $m$ are related as
$$(3.4) \qquad i = i_m, \qquad j = \mu_{i,m}.$$

Clearly, whichever choice of $(I_1, I_2, \ldots, I_M)$ is given, for any $(i, j)$ with $j \leq \nu_i$, there is exactly one $m$, $1 \leq m \leq L$, that satisfies (3.4) and vice versa. Note that here $m$ is a random variable (it depends on $I_1, \ldots, I_M$).

Since $\nu_i \geq \mu_{i,m}$, we have
$$(3.5) \qquad E\{\zeta_{ij}^2 | \nu_i \geq j\} = E\{\eta_m^2 | \nu_i \geq \mu_{i,m}\} = E\eta_m^2.$$

We have
$$\sum_m \eta_m^2 = \sum_{i,j} \zeta_{ij}^2 = \sum_i \sum_{i \in I_l} |Z_i(l) - Z_i(l-1)|^2$$
$$= \sum_{l=1}^M |Z_{I_l}(l) - Z_{I_l}(l-1)|^2.$$

This, together with (3.3) and (3.5), completes the proof of Lemma 1. □

PROOF OF THE AUXILIARY THEOREM. By Lemma 1, all we have to prove is that
$$E \sum_{l=1}^M |Z_{I_l}(l) - Z_{I_l}(l-1)|^2 \leq \frac{2}{\rho} \cdot D(p^n \| q^n).$$

In fact, we prove that, for any realization

(3.6) $$I_1, I_2, \ldots, I_M$$

of the sequence of patches, we have

(3.7) $$E\left\{ \sum_{l=1}^M |Z_{I_l}(l) - Z_{I_l}(l-1)|^2 \Big| I_1, I_2, \ldots, I_M \right\} \leq \frac{2}{\rho} \cdot D(p^n \| q^n).$$

The left-hand side of (3.7) can be written as
$$\sum_{l=1}^M E\{|Z_{I_l}(l) - Z_{I_l}(l-1)|^2 | I_1, I_2, \ldots, I_l\}.$$

Fix the sequence (3.6) and recall from Section 2 the definition
$$r^n(l) = \operatorname{dist}(Z^n(l) | I_1, I_2, \ldots, I_l),$$
according to which $r^n(l)$ is obtained from $r^n(l-1)$ by putting

(3.8) $$\operatorname{dist}(Z_{I_l}(l) | \bar{z}_{I_l}) = Q_{I_l}(\cdot | \bar{z}_{I_l})$$

and leaving unchanged the distribution of the coordinates outside $I_l$:

(3.9) $$(\bar{r}(l))_{I_l} = (\bar{r}(l-1))_{I_l}.$$

It follows from (2.3), (3.8) and (3.9) that the joint conditional distribution
$$\operatorname{dist}(Z^n(l-1), Z^n(l) | I_1, \ldots, I_l)$$
can be defined in such a way that

(3.10) $$E\{|Z_{I_l}(l) - Z_{I_l}(l-1)|^2 | I_1, \ldots, I_l\} \leq \frac{2}{\rho} \cdot D(r^n(l-1) \| r^n(l)).$$

[The distributions $r^n(l-1)$ and $r^n(l)$ depend on $I_1,\ldots,I_l$.] Indeed, the left-hand side of (3.10) can be written as

$$\int E\{|Z_{I_l}(l) - Z_{I_l}(l-1)|^2|\bar{z}_{I_l}(l-1), I_1, I_2, \ldots, I_l\}(\bar{r}(l-1))_{I_l}(\bar{z}_{I_l})\,d\bar{z}_{I_l}$$

$$\leq \frac{2}{\rho} \cdot \int D(r_{I_l}(l-1)(\cdot|\bar{z}_{I_l}(l-1))\|Q_{I_l}(\cdot|\bar{z}_{I_l}(l-1)))(\bar{r}(l-1))_{I_l}(\bar{z}_{I_l})\,d\bar{z}_{I_l}$$

$$= \frac{2}{\rho} \cdot D(r^n(l-1)\|r^n(l)).$$

The last equality here follows from (3.9).

Therefore, it is enough to prove that for any choice of $I_1,\ldots,I_M$,

$$(3.11) \qquad D(p^n\|q^n) \geq \sum_{l=1}^{M} D(r^n(l-1)\|r^n(l)).$$

This follows from the identities

$$(3.12) \quad \begin{aligned} D(p^n\|q^n) &= D(p^n\|r^n(1)) + D(r^n(1)\|r^n(2)) + \cdots \\ &\quad + D(r^n(l-1)\|r^n(l)) + D(r^n(l)\|q^n), \end{aligned}$$

valid for any $l \geq 1$. It is clear that (3.12) for $l = M$ implies (3.11).

We prove (3.12) by induction on $l$. Thus first we claim that

$$(3.13) \qquad D(p^n\|q^n) = D(p^n\|r^n(1)) + D(r^n(1)\|q^n),$$

which is just (3.12) for $l=1$. Indeed, by the well-known decomposition formula for divergence,

$$D(p^n\|q^n) = D(\bar{p}_{I_1}\|\bar{q}_{I_1}) + \int \log \frac{p_{I_1}(y_{I_1}|\bar{y}_{I_1})}{Q_{I_1}(y_{I_1}|\bar{y}_{I_1})} p^n(y^n)\,dy^n$$

$$= D(r^n(1)\|q^n) + D(p^n\|r^n(1)).$$

Now apply (3.13) to $r^n(1)$ in the role of $p^n$. This, together with (3.13), yields

$$D(p^n\|q^n) = D(p^n\|r^n(1)) + D(r^n(1)\|r^n(2)) + D(r^n(2)\|q^n).$$

Iterating this step, (3.12) follows for any $l$. □

**4. Proof of Proposition 2.** Consider the joint distribution $\text{dist}(Y^n, U^n)$, achieving $W^2(p^n, r^n)$. Let $Y^n(1)$ and $U^n(1)$ denote random variables with density functions $p^n\Gamma$ and $r^n\Gamma$, respectively.

For a given $I \in \mathcal{I}$, we define a joint conditional density function

$$\text{dist}(Y^n, U^n, Y^n(1), U^n(1)|I)$$

as follows. Put $\bar{Y}_I(1) = \bar{Y}_I$, $\bar{U}_I(1) = \bar{U}_I$,

(4.1) $$\text{dist}(Y_I(1)|Y^n = y^n, I) = Q_I(\cdot|\bar{y}_I),$$

(4.2) $$\text{dist}(U_I(1)|U^n = u^n, I) = Q_I(\cdot|\bar{u}_I),$$

and take for
$$\text{dist}(Y_I(1), U_I(1)|Y^n = y^n, U^n = u^n, I)$$
a joining of (4.1) and (4.2) to achieve
$$E\{|Y_I(1) - U_I(1)|^2|Y^n = y^n, U^n = u^n, I\} = W^2(Q_I(\cdot|\bar{y}_I), Q_I(\cdot|\bar{u}_I)).$$

We have by (1.3)
$$E|Y^n(1) - U^n(1)|^2$$
$$\leq \frac{1}{N} \sum_{I \in \mathcal{I}} E\left[\sum_{k \notin I} |Y_k - U_k|^2 + W^2(Q_I(\cdot|\bar{Y}_I), Q_I(\cdot|\bar{U}_I))\right]$$
$$\leq \frac{1}{N} \sum_{I \in \mathcal{I}} \sum_{k \notin I} E|Y_k - U_k|^2 + \frac{1}{N}(1-\delta)tE|Y^n - U^n|^2$$
$$\leq \left(1 - \frac{t}{N}\right) \sum_{k=1}^{n} E|Y_k - U_k|^2 + \frac{1}{N}(1-\delta)tE|Y^n - U^n|^2$$
$$= \left(1 - \frac{t\delta}{N}\right) E|Y^n - U^n|^2.$$

Proposition 2 is proved.

**5. Proof of Theorem 1.** Let $M$ be fixed, and apply Proposition 2 $M$ times to the distributions
$$p^n = \text{dist}\, Y^n \quad \text{and} \quad r^n = p^n \Gamma^M = \text{dist}\, Z^n.$$
(We use the notation of Section 2.) We get that
$$W^2(p^n \Gamma^M, r^n \Gamma^M) = W^2(p^n \Gamma^M, p^n \Gamma^{2M})$$
$$\leq \left(1 - \frac{t\delta}{N}\right)^M \cdot W^2(p^n, r^n) \leq \exp\left(-t\delta\frac{M}{N}\right) \cdot W^2(p^n, r^n),$$
that is,
$$W(p^n \Gamma^M, r^n \Gamma^M) = W(p^n \Gamma^M, p^n \Gamma^{2M}) \leq \exp\left(-t\delta\frac{M}{2N}\right) \cdot W(p^n, r^n).$$

Iterating this step, we get that for any $j \geq 1$,
$$W(p^n \Gamma^{(j-1)M}, p^n \Gamma^{jM}) \leq \exp\left(-(j-1)t\delta\frac{M}{2N}\right) \cdot W(p^n, r^n).$$

Let us define the random sequences $Y^n(j)$, $j = 0, 1, \ldots$, so that

$$\operatorname{dist} Y^n(j) = p^n \Gamma^{jM}$$

and, for $j \geq 1$,

(5.1) $\quad [E|Y^n(j) - Y^n(j-1)|^2]^{1/2} \leq \exp\left(-(j-1)t\delta\frac{M}{2N}\right) \cdot W(p^n, r^n).$

We see that $\{Y^n(j)\}$ is a Cauchy sequence in $L_2$ and thus it converges in $L_2$ to some random sequence $\tilde{X}^n$. However, we must have $\operatorname{dist} \tilde{X}^n = q^n$. Indeed, $q^n$ is invariant with respect to $\Gamma$ and therefore Proposition 3 implies that $W(p^n \Gamma^{jM}, q^n) \to 0$ as $j \to \infty$. Thus we can assume that the sequence $\{Y^n(j)\}$ converges to $X^n$ in $L_2$.

By the estimates (5.1),

$$[E|Y^n - X^n|^2]^{1/2} \leq [E|Y^n - Z^n|^2]^{1/2} \cdot \frac{1}{1 - \exp(-t\delta(M/(2N)))}.$$

By the Auxiliary Theorem, this implies

$$[E|Y^n - X^n|^2]^{1/2} \leq \sqrt{\frac{M}{N} \cdot v \cdot \frac{2}{\rho} \cdot D(p^n \| q^n)} \cdot \frac{1}{1 - \exp(-t\delta(M/(2N)))}$$

$$= \sqrt{2} \cdot \frac{\sqrt{t\delta(M/(2N))}}{1 - \exp(-t\delta(M/(2N)))} \cdot \sqrt{\frac{1}{\delta} \cdot \frac{v}{t} \cdot \frac{2}{\rho} \cdot D(p^n \| q^n)}.$$

Now to complete the proof, it is enough to see that the factor

$$\frac{\sqrt{t\delta(M/(2N))}}{1 - \exp(-t\delta(M/(2N)))}$$

can be bounded by a numerical constant through an appropriate selection of $M$. Notice that the function

$$f(x) = \frac{\sqrt{x}}{1 - e^{-x}}, \qquad x > 0,$$

is bounded in any bounded interval that is bounded away from 0. If $M$ varies on the integers, then the quantity $x = t\delta M/(2N)$ changes by steps smaller than $1/2$. Thus there is a value of $M$ for which $x = t\delta M/(2N)$ is between 1 and $3/2$, and so

$$\min_M f(t\delta M/(2N)) \leq \max_{1 \leq x \leq 3/2} f(x).$$

This completes the proof of Theorem 1.

**6. Proof of Theorem 2.** The Otto–Villani theorem implies condition (1.5) of Theorem 1. To prove the contractivity condition, fix two sequences $x^n$, $y^n \in \mathbb{R}^n$. By (1.5) and the logarithmic Sobolev inequality, we have

$$W^2(Q_I(\cdot|\bar{x}_I), Q_I(\cdot|\bar{y}_I)) \leq \frac{1}{2\rho} \cdot D(Q_I(\cdot|\bar{x}_I) \| Q_I(\cdot|\bar{y}_I))$$

$$\leq \frac{1}{\rho^2} \cdot \int_{\mathbb{R}^I} \sum_{i \in I} [\partial_i \Phi(\eta_I \bar{x}_I) - \partial_i \Phi(\eta_I \bar{y}_I)]^2 Q(\eta_I|\bar{x}_I) \, d\eta_I.$$

It follows that

$$\sum_{I \in \mathcal{I}} W^2(Q_I(\cdot|\bar{x}_I), Q_I(\cdot|\bar{y}_I))$$

(6.1)
$$\leq \frac{1}{\rho^2} \cdot \sum_{I \in \mathcal{I}} \int_{\mathbb{R}^I} \sum_{i \in I} [\partial_i \Phi(\eta_I \bar{x}_I) - \partial_i \Phi(\eta_I \bar{y}_I)]^2 Q(\eta_I|\bar{x}_I) \, d\eta_I$$

$$= \frac{1}{\rho^2} \cdot \int \sum_{I \in \mathcal{I}} \sum_{i \in I} [\partial_i \Phi(\eta_I \bar{x}_I) - \partial_i \Phi(\eta_I \bar{y}_I)]^2 \prod_{I \in \mathcal{I}} Q_I(\eta_I|\bar{x}_I) \prod_{I \in \mathcal{I}} d\eta_I.$$

(The integral in the last line is taken over $\prod_{I \in \mathcal{I}} \mathbb{R}^I$.)

Now consider, for a fixed vector $\eta = (\eta_I, I \in \mathcal{I})$, the function

$$g = g^\eta : \mathbb{R}^n \mapsto \prod_{I \in \mathcal{I}} \mathbb{R}^I,$$

defined by

$$g_{I,i}(y^n) = \frac{1}{\rho} \cdot \partial_i \Phi(\eta_I \bar{y}_I), \qquad i \in I \in \mathcal{I}.$$

Observe that the expression

$$\frac{1}{\rho^2} \sum_{I \in \mathcal{I}} \sum_{i \in I} [\partial_i \Phi(\eta_I \bar{x}_I) - \partial_i \Phi(\eta_I \bar{y}_I)]^2,$$

integrated (with respect to some density function) in the last line of (6.1), is nothing else than the squared Euclidean norm of the increment of $g^\eta$ between the points $x^n$ and $y^n$. By assumption (1.9) of Theorem 2, the norm of the Jacobian of $g^\eta$ is bounded by $(t \cdot (1-\delta))^{1/2}$, so (6.1) implies the contractivity condition (1.3).

**Acknowledgment.** I thank the Erwin Schrödinger Institute in Vienna for its hospitality and M. Ledoux, who called my attention to this problem during my stay at the Erwin Schrödinger Institute in 1999.

Alfréd Rényi Institute of Mathematics
Hungarian Academy of Sciences
H-1364 Budapest
P.O. Box 127
Hungary
e-mail: marton@renyi.hu
url: www.renyi.hu/~marton